\documentclass{jgcc} %%% last changed 2020-01-17


\usepackage{hyperref}
\usepackage{amssymb}
\theoremstyle{plain} %\crefname{satz}{Satz}{S\"atze}
\usepackage{graphicx,multirow,diagbox,url,booktabs,algpseudocode,algorithm}
\newcolumntype{M}[1]{>{\centering\arraybackslash}m{#1}}%%%%%%%

\begin{document}

\title[Hyper-heuristics in Group Cryptology]{Evolution of Group-Theoretic Cryptology Attacks using Hyper-heuristics}

\author[M.~J.~Craven]{Matthew J. Craven}
\address{Centre for Math. Sciences, University of Plymouth, Drake Circus, Plymouth, PL4 8AA, U. K.}
\email{matthew.craven@plymouth.ac.uk}

\author[J.~R.~Woodward]{John R. Woodward}
\address{Queen Mary University of London, Mile End Rd, London E1 4NS, U. K.}
\email{j.woodward@qmul.ac.uk}

%\classification{20P05, 68W30, 90C27, 94A60}

\begin{abstract}
  \noindent In previous work, we developed a 
\emph{single} Evolutionary Algorithm (EA) to solve random instances of the Anshel-Anshel-Goldfeld (AAG) key exchange protocol over polycyclic groups. 
The EA consisted of six simple heuristics which manipulated strings. 
The present work extends this by
exploring the use of hyper-heuristics in group-theoretic cryptology for the first time. 
Hyper-heuristics are a way to generate new algorithms from existing algorithm components
(in this case the simple heuristics),
with the EAs being one example of the type of algorithm
which can be generated by our hyper-heuristic framework. 
We take as a starting point the above EA
and allow hyper-heuristics to
build on it by making small tweaks to it. 
This adaptation is through a process of taking the EA and injecting chains of heuristics 
built from the simple heuristics. 

We demonstrate we can create novel heuristic chains, 
which when placed in the EA create algorithms which out-perform the existing EA. 
The new algorithms solve a markedly greater number of random AAG instances than the EA for harder instances. 
This suggests the approach could be applied to many of the same kinds of problems, 
providing a framework for the solution of cryptology problems over groups.
The contribution of this paper is thus a framework 
to automatically build algorithms to attack cryptology problems.
\end{abstract}

\maketitle

\section{Introduction}\label{sec:intro}
On NP-hard problems, 
the time taken to produce an algorithm to solve such problems is often vast. 
In such cases, users may use an ``off the shelf'' algorithm to 
obtain approximate solutions within an appropriate time.
In this paper, we take a different
approach and attempt to design
an algorithm in response to feedback from 
similar instances of the problem. 
Examples of such problems are those in group-theoretic cryptology 
(multiple conjugacy, Anshel-Anshel-Goldfeld (AAG, \cite{anshel2001new}) and word decomposition, for instance).
These problems have been posed over varying types of groups serving as the base problems for key exchange protocols (KEPs) \cite{anshel1999algebraic,anshel2001new,ko2000new,eick2004polycyclic,kahrobaei2014heisenberg} and subsequently attacked \cite{garber2005probabilistic,garber2006length,franco2003conjugacy,craven2012evolutionary,craven2016parallel,kotov2015analysis,myasnikov2007length,myasnikov2008group,ruinskiy2007length}. 
The group structures used are often intended to provide an extra encryption layer through the scrambling induced by the group presentation.

In this work,
a preliminary hyper-heuristic framework is detailed which takes as input a proposed cryptographic base problem 
and a group structure, and, 
via machine learning techniques, generates operations for a length attack algorithm which aims to solve an acceptable proportion of random instances of the problem. 
The framework is implemented in the GAP 4.8.7 \cite{GAP4} language 
(due to its compatibility with the ParGAP package \cite{ParGAP}, 
allowing use of MPI intra-core communication).
This is tested on a case study of an AAG KEP 
\cite{anshel1999algebraic,anshel2001new} posed over polycyclic groups defined by a number field \cite{garber2015length}. 
The aim is to generate mutation operators for algorithms which outperform the existing human-designed EA. 
These mutation operators are chains of simple heuristics,
which are composed or learned. 
The generation of crossover, selection, and other heuristic components are outside the scope of this work. 

Our contribution is an approach that, contrary to the above manual design of attacks, automatically builds attack mechanisms and attempts to break the above AAG KEP. 
This approach is trained on small set of instances and then validated on a second larger independent set of instances, illustrating it generalises.
This paper is not proposing a single algorithm to attack, but rather a framework
in which algorithms can automatically be generated and then tested, 
and is an example of the generate-and-test paradigm which has many applications in science, engineering, mathematics, and daily life. 
One of the drawbacks of our method is the large amount of computation time required; it takes relatively little time to generate an algorithm but a relatively long time to test it. 
One of the benefits of our approach, however, is that we can take existing approaches (as we do in this paper, an EA 
\cite{woodward2011automatically, woodward2012automatic}), and use it as a starting point from which we can improve. 

In \cite{craven2016parallel}
we observed that a human-designed EA performs better than the length attack algorithm of \cite{garber2015length}. 
In this paper, we observe that an automatically designed EA performs better than the human designed EA.  
We also conjecture that a random search algorithm will perform poorly on this problem. 
This is a pattern of performance typically seen in the metaheuristics literature. 
The reason for this ordering of four types of solver lies in the nature of the resulting search landscape. 
A human designed EA is essentially a more sophisticated length attack algorithm, 
and a machine-designed EA is essentially slightly more sophisticated than an human designed EA. 

Typically, during the design of an algorithm, we need an understanding of the problem to design an algorithm.
The algorithm thus capture our intuition about how to solve that problem.
(consider the problem of sorting and the large number of algorithms available, for instance).
An algorithm is an explicit formalisation of our intuition: with cryptology, we have very little in the way of intuition to guide us.
This is an opportunity for an automated method (which is largely unbiased) to invent new algorithms.

It is acknowledged that the detailed protocol has already been broken by \cite{craven2016parallel,kotov2015analysis}
(the latter reference being a ``field attack''),
but the wish is to present this work as a preliminary study with a view towards 
application to other cryptanalytic problems. 
It is argued that this type of algorithm has a future in the disciplines of cryptology and possibly algorithmic questions in combinatorial group theory, and may be extended to other structures and problem types.

This work is organised as follows:
In Section \ref{Background} we give an introduction to group-based cryptography, reviewing previously-proposed KEP problems, before turning to an overview of hyper-heuristics.
This is followed by Section \ref{Notation and Formalisation}, which 
introduces the notation and formalisation.
In Section \ref{Experimental Setup} we describe the experimental approach and detail parameter settings, discussing the results of our approach in Section \ref{experimental results}.
In Section \ref{Concluding remarks} we conclude the article, including a discussion of further work resulting from this study and raising future research directions.  

\section{Background}
\label{Background}
In this section we will first introduce group-based cryptography. We then give an introduction to hyper-heuristics. 

\subsection{Introduction to Group-Based Cryptography}
Group-based cryptography uses groups in the construction of cryptosystems and KEPs and has been an active area of research since approximately the late 1990s. Proposed cryptosystems and their subsequent attacks (purported breaks) iterate one after the other with the aim of producing increasingly secure cryptography over time.

The late nineties were when group-based cryptography began in earnest, when the likes of \cite{anshel1999algebraic,anshel2001new,ko2000new} proposed KEPs based upon braid groups. As mentioned in the introduction, the braid groups were used due to the scrambling induced by the presentation of the group, and the consequent thought that the underlying problems (various guises of the conjugacy problem) were thought to be extremely difficult to solve. Solving the underlying problem would, in many cases, break the KEP and render any keys exchanged open to misuse by adversaries.

Both KEPs, and the underlying problems, were attacked in the next few years. 
Example of such attacks were super summit set attacks \cite{franco2003conjugacy} and the more practical length-based attacks (LBAs) \cite{hughes2003length}.
These latter algorithms (also known as hillclimbers) build up solutions to instances of the problem gradually, beginning with a short candidate solution and making alterations to it based upon randomness. 
This altered solution is then compared to the old solution by some metric, mostly with regards to how ``well'' the candidate solves the instance 
(for example, how many symbols remain after all possible cancellations have been conducted). 
If the altered solution proves to be an improvement then the current solution is set equal to the altered solution and the process is repeated. 
If not, then the altered solution is discarded.

Being practical and fast, LBAs became increasingly sophisticated through \cite{garber2006length,garber2005probabilistic,myasnikov2007length,ruinskiy2007length}. As LBAs became also increasingly capable of solving instances of the aforementioned KEPs, 
researchers began, in a search for more attack-resistant structures, to look for new groups and problems while keeping the general methodology.
Examples of these platform groups are right-angled Artin groups \cite{craven2012evolutionary} (a homomorphic pre-image of braid groups), small cancellation groups \cite{shpilrain2009using}, matrix groups, Thompson's group and Grigorchuk's group, to name but a few.

Polycyclic groups were first proposed as a new platform group in 2004 \cite{eick2004polycyclic} and were followed ten years later by the works of \cite{kahrobaei2014heisenberg} and \cite{garber2015length}, applying two distinct types of polycyclic groups to the AAG \cite{anshel2001new} problem (multiple conjugacy). 
The systems introduced were, in turn, broken by the works of \cite{blaney2016ptime} (for generalised Heisenberg groups),
\cite{kotov2015analysis} and \cite{craven2016parallel} (via a parallelised EA). 
The latter work was demonstrated to be more efficient, and more successful, than previous LBA attacks. 
Although the approach on the proposed KEP was successful, we wish to take it further into the domain of hyper-heuristics and use the KEP as a test bed for our framework. An excellent summary of group-theoretic cryptology in general can be found in \cite{myasnikov2008group}.

\subsection{Introduction to Hyper-Heuristics} 
Informally, hyper-heuristics offer to take a number of existing computational search techniques, and combine them, to make a new heuristic. 
This new heuristic is intended to have more of the strengths of each of the heuristics, and less of their weaknesses.
The motive of a hyper-heuristic is not to out-perform a state-of-the-art algorithm on a single instance of a problem.
Rather, the aim of hyper-heuristic approaches is to perform well across a range of problem instances. 
In other words, hyper-heuristics attempt to offer robust performance across a set of problems rather than specialised performance on a narrow set of specific instances. 
These problems could be problem instances from a given domain, such as the travelling salesman problem. 
Or the problem instances could be drawn from different problem domains, for example exam timetabling and vehicle routing. In this paper we are developing a hyper-heuristic framework to solve problem instances from a single domain: cryptology.
 
We should also be careful about the distinction between optimisation and supervised machine learning.
Optimisation typically has an objective function we wish to evaluate and a parameter value which is a global optimum. 
Often this is difficult to achieve, and also difficult to know when it has been achieved. 
In contrast, with supervised machine learning, we typically have a set of example cases which we use to train a model.
We then have a second set of independent example cases which are used to determine if the model performs well in general on cases which were not included in the training phase.
Optimisation has a single stage (optimising), while machine learning has two main stages (training and testing).
Nor do we have the issue of over-fitting in optimisation, but the issue of over-fitting may arise in machine learning. 
In summary, in this paper we are using a machine learning approach (hyper-heuristics), with an independent training and test set, to build a heuristic which we used for optimisation, the objective function being to minimise the length.

Hyper-heuristics can be viewed in the context of heuristics and metaheuristics. These three terms are often confused. Let us begin by looking first at heuristics, metaheuristics, and finally hyper-heuristics.

A heuristic is domain-specific algorithm (often called a rule of thumb) which does not solve a problem to optimality
(as such problems are often NP-hard or NP-complete), but rather offers to deliver suboptimal solutions in feasible time. 
That is, a heuristic is a strategy that aims to deliver an approximation to a solution to a given problem in a fast,
rather than an overly elaborate, way. 
An example of a heuristic is the Lin-Kernighan algorithm which is applied to the Traveling Salesman Problems (TSP). 
It does not make sense to apply the Lin-Kernighan algorithm to the knapsack problem, as it is specific to TSP problems.
The Lin-Kernighan algorithm could be applied to other graph-based problems with a representation similar to the TSP, but the algorithm may not perform well as this is not what it was intended for. A metaheuristic is a general search-based algorithm which can be applied to spaces consisting of bit strings or permutations, for example, depending on the representation of the problem instances. An example of a metaheuristic is a genetic algorithm which searches the space of bit strings of a given length.

Hyper-heuristics are different again. Typically a hyper-heuristic uses a metaheuristic to search the space of problem specific heuristics. 
That is, a hyper-heuristic is a ``search methodolog[y] for choosing or generating (combining, adapting) heuristics [...], in order to solve a range of optimisation problems'' \cite[p.~2]{burke2009exploring}. 
For example, see \cite{Bai07amodel}.
Hyper-heuristics have successfully been applied to a number of different problem domains.

As combinatorial optimisation problems are a subset of all NP hard problems, 
it is not surprising that hyper-heuristics have been a popular approach. 
Applications include 
exam timetabling \cite{bilgin2006experimental}
bin packing 
\cite{ross2002hyper}
and employee rostering
\cite{burke2003tabu}.
There have also been a number of well-referenced survey articles, including 
\cite{ross2005hyper, burke2013hyper, burke2019classification}.

Hyper-heuristics typically do not generate complete algorithms; 
rather a component of an algorithm is targeted to be automatically designed
by a generate and test approach. Hyper-heuristics have been used, for example, to generate
components of evolutionary algorithms such as 
genetic algorithms and evolutionary programming
(e.g., 
crossover operators \cite{goldman2011self}, mutation operators \cite{hong2013automated}) and form a large part of the literature in the automated design of algorithms \cite{haraldsson2014automated}.

In the context of this paper, we are using hyper-heuristics in the following manner.
We will take seven low-level heuristics, which are chained together randomly to effectively create new heuristics.
These new chains of heuristics are then inserted into a standard EA (depicted in the work of \cite{craven2016parallel}) which is used to tackle the problem. This work begins in the next section.

\section{Notation and Formalisation}
\label{Notation and Formalisation}

In this section, the AAG KEP over a certain type of polycyclic group is discussed. 
This is followed by the notation needed for the implementation of the hyper-heuristic. 
In this section, the notation broadly follows that of \cite{craven2016parallel} which describes the aforementioned EA. 

\subsection{Setup of Problem}
\label{defn_problem}
The AAG KEP \cite{anshel1999algebraic,anshel2001new} was posed over polycyclic groups in \cite{garber2015length,kahrobaei2014heisenberg},
and subsequently attacked in two distinct ways by the work of \cite{craven2016parallel} and \cite{kotov2015analysis}. 
The main details of the protocol, following the exposition given in \cite{craven2016parallel} for a group $G=\left<g_1,g_2,\ldots,g_n\,|\,R\right>$, are as follows.

First, 
Alice chooses a subgroup $\mathcal{A}=\left<a_1,a_2,\ldots,a_N\right>\leq G$ generated by words $a_i$ in the generators of $G$ such that $L_1\leq l_G(a_i)\leq L_2$. 
Bob then does similarly to produce a subgroup $\mathcal{B}=\left<b_1,b_2,\ldots,d_N\right>\leq G$. 
All of $A$, $B$ and $G$ are made public. 
Alice chooses her private key $A=a_{\mu_1}^{\epsilon_1}a_{\mu_2}^{\epsilon_2}\ldots a_{\mu_L}^{\epsilon_L}$ where each $a_i\in \mathcal{A}$ and $\mu_i=\pm 1$ (for all $i=1,\ldots,L$). 
She now calculates
\[
A^{-1}b_1 A, A^{-1}b_2 A, \ldots, A^{-1}b_N A,
\]
and sends these to Bob. 
Bob does similarly, producing $B^{-1}a_i B$ for $i=1,\ldots,N$ and sends these to Alice (his private key is $B$). 
From the information now exchanged, each individual can now produce the shared key (the commutator) $[A,B]=A^{-1}B^{-1}AB$.

If an adversary wishes to find either the private key $A$ (or equivalently, $B$), they may intercept the above conjugates either party sends to the other. 
Thus the 
problem to be solved may be simply expressed as a subgroup restricted multiple conjugacy problem in the following way. 
Each instance of this problem is a set of $N$ (frequently twenty) conjugacy equations $E = \{E_1, \ldots, E_N\}$
\begin{eqnarray*}
E_1: c_1 &=& A^{-1}b_iA\\
E_2: c_2 &=& A^{-1}b_iA\\
&\ldots&\\
E_N: c_N &=& A^{-1}b_NA,
\end{eqnarray*}
posed over a finitely presented platform group $G$. 
A solution to the problem means that all the above equations are satisfied.
One function of the rewriting rules (relators) $R$ of $G$ is to serve cryptographically as word obfuscators and thus hide the secret word (private key) $A$.

The problem is posed over polycyclic groups $\mathcal{O}\rtimes U$, where, by \cite{craven2016parallel}, $\mathcal{O}$ is the additive group of the ring of integers of a number field $K$ and its group of units is $U$. 
The number field is written as $K=\mathbb{Q}[x]/(f)$, for $f\in\mathbb{Z}[x]$ a monic irreducible polynomial of degree $d$.
To recap, the instance parameters associated to this setup are then the number of equations, $N$, the polynomial $f$, length $L$ of the private key $A$ in $\mathcal{A}$, and $L_1$ and $L_2$ (the lower and upper bounds, respectively, on the lengths of $a_i$ in $G$).

Note that, in this work, we refer to either an \emph{exact} solution or a \emph{candidate} solution as appropriate.
However, most references will be to candidate solutions but for the sake of brevity will be named \emph{solutions}.
In the context of hyper-heuristics and cryptology,
we are using hyper-heuristics to generate candidate solutions to find an exact solution to the cryptographic problem. In this spirit, there are several functions at work in this paper which we need to distinguish between.

\subsection{Pertinent Functions}
The following functions are recapped from \cite{craven2016parallel}.
Let a word $w$ be expressed in the form $w=f_{i_1}^{e_{i_1}}f_{i_2}^{e_{i_2}}\ldots f_{i_r}^{e_{i_r}}$ 
for non-zero $e_j\in\mathbb{Z}$ and $f_1$, $f_2$, $\ldots$, $f_n$ are the generators of the free group $F$.
The \textit{length functions} associated to the group $G$ are then given by
\[
\ell\left(w\right)=\sum\limits_{k=1}^r |e_{i_k}|
\]
and
\[
\ell_{wt}\left(w\right)=\sum\limits_{k=1}^r \omega_{i_k} |e_{i_k}|,
\]
where, as in the above, $\omega_j$ is the ``sum of the lengths of the normal forms of the commutators $[g_j,g_k]$ in $G$ for $k=1,\ldots,n$''. 
That is, the length of $w$ is the sum of the absolute powers (respectively, the weighted absolute powers) of individual generators $f_i$ that make up the word $w$.

The basic EA \textit{cost function} measures the quality of the candidate solutions produced by the EA and is given by
\[
c\left(w\right)=\sum\limits_{i=1}^N \ell\left(\alpha^{-1}b_i\alpha c_i^{-1}\right),
\]
where $\alpha$ is the current EA solution (i.e., the approximation of the private key $A$).
This function has output of the sum of lengths of (normal form) reduced equations $E_1,E_2,\ldots,E_N$. 
That is, the length of summand $i$ (where $i\in\left\{1,\ldots,N\right\}$) of the cost function is equal to the reduced length of each equation $E_i$ after its substitution with $\alpha$. 
This function is used to drive search in the EA, since the population ranking is performed with respect to it.
The cost used in the EA is broadly the cost vector produced by this basic function, involving the sum $c$, maximum and mean lengths of summands of $c$ for the weighted and non-weighted length functions, given in \cite[p.~8--9]{craven2016parallel}. 
The global optimum (minimum value) of $c$ is zero; at this value, no fragments of the equations remain and the instance is completely solved.

The \textit{heuristic objective function} is the metric used to compare the current heuristic chain over the given set of instances (training or testing) and is a vector given in the following order, each element computed over the set of instances:
\begin{itemize}
\item The mean best cost $c$ over the unsuccessful EA runs;
\item The negative of the success rate as a proportion of the total number of runs;
\item From the successful runs, the mean number of generations used by the EA.
\end{itemize}
That is, this function tells the hyper-heuristic how good a given heuristic chain is. 
For the validation process the first and second elements of the above objective function are swapped, since we are now more concerned with the success rate. 
The hyper-heuristic attempts to minimise the above objective function, indicating a successful heuristic chain, as far as possible. 
Comparison of heuristic objective vectors, produced by two distinct heuristic chains, is performed lexicographically. Note that this function is often termed a fitness function in the evolutionary computation community.

\subsection{Simple Heuristics on the Group}\label{sec:simple Heuristics on the group}

In previous work \cite{craven2016parallel},
six simple heuristics were used in an EA to 
break a proposed key exchange \cite{garber2015length}. 
These are listed in Table \ref{simple_heuristics} as $H_1$--$H_6$. 
In this paper, we are building 
new heuristic chains to 
inject into an EA. 
We have also added a seventh heuristic $H_7$ (swap) to this set of heuristics. 
Evolutionary operators may be otherwise thought of as heuristic on group elements $w=f_{i_1}^{e_{i_1}}f_{i_2}^{e_{i_2}}\ldots f_{i_r}^{e_{i_r}}$.

\begin{table}[h!]
\begin{center}
\begin{tabular}{cc}
\hline
Heuristic & Description \\
\hline
$H_1$ & insertion of a subgroup generator $a_i$ \\
$H_2$ & insertion of a single generator $f_i$ \\
$H_3$ & deletion of a single generator \\
$H_4$ & substitution with a single generator \\
$H_5$ & position conjugation: conjugating a given position by $f_i$ \\
$H_6$ & subword conjugation: conjugating a subword by $f_i$ \\
$H_7$ & swap \\
\hline
\end{tabular}
\caption{The seven simple heuristics used to build new heuristics. The first six heuristics were used in \cite{craven2016parallel}.}\label{simple_heuristics}
\end{center}
\end{table}

Heuristic $H_7$ (swap) is designed to assist when symbols are in the `wrong place' in a word $w$, swapping two symbols at random positions and potentially trigger subsequent cancellation of symbols (and, thus, an EA cost reduction).
Essentially all heuristics in the above table are random, with operations performed with random words or generators at random positions.
The above is not a list of minimal heuristics: it is noted, for example, that heuristic $H_1$ can be achieved through repeated application of $H_2$, as can $H_5$ and $H_6$ (which were specialised to the conjugacy problem).

\subsection{EA Parameter Settings}
The EA parameters are given in Table \ref{new_EA_params}. 
The parameters were produced by copious experimentation, and scaling down the parameters in \cite{craven2016parallel} to approximately one quarter of their original values to achieve an effective set of EA parameters. 
This increases the speed of the EA. 
The original population size was 100. We do not claim optimality for these parameter settings. 

\begin{table}[h!]
\begin{center}
\begin{tabular}{cc}
\hline
Parameter & Value \\
\hline
population size & 25 \\ 
truncation selection & 40\% \\  
$H_1$ & 6 \\
$H_2$ & 1 \\
$H_3$ & 1 \\
$H_4$ & 5 \\
$H_5$ & 1 \\
$H_6$ & 1 \\ 
crossover & 4 \\  
selection & 2 \\
chains & 4 \\
number of generations & depends on experiment\\
\hline
\end{tabular}
\caption{EA parameter settings.}\label{new_EA_params}
\end{center}
\end{table}

All heuristics are performed by firstly choosing a solution at random from the top 40\% of the population (by cost). 
The selection operator is elitist; i.e., if $n_s$ solutions are to be selected then the first is the solution of minimum cost with the remaining $n_s-1$ solutions selected from the top 40\% (i.e., after ranking by minimum cost) of the population at random. 
All random choices are made uniformly (as in \cite{craven2016parallel}).

It was chosen to have four solutions from each generation created by a heuristic chain.
Testing this alongside the remaining nineteen solutions in each generation created by the same heuristic chain, it was found that this choice of four solutions turned out to be more advantageous (the average number of generations to solve decreased). 
$H_7$ does not appear in Table \ref{new_EA_params}
as it does not operate in isolation (as part of the EA of \cite{craven2016parallel}), 
only in the context of the other six heuristics. 
Crossover is performed by choosing two words (from the top 40\% of the population) $w_1$, $w_2$. 
Choosing two random positive integers $r_1\leq \ell(w_1)$, $r_2\leq \ell(w_2)$, one of the two words
\[
w_1[1\ldots r_1]w_2[r_2+1 \ldots \ell(w_2)]\,\,\,\,\,\textrm{and}
\,\,\,\,\,w_2[1\ldots r_2]w_1[r_1+1\ldots \ell(w_1)]
\]
is output \cite{craven2016parallel} (where $w[s\ldots t]$ is the subword between and including positions $s$ and $t$ of the word $w$). 
The next section details the operation of the hyper-heuristic and the experimental setup.

\section{Experimental Setup}
\label{Experimental Setup}

\subsection{Hyper-heuristic Implementation}
As above, our objective is to create a hyper-heuristic that, 
given the AAG problem (Section \ref{defn_problem}) and a polycyclic group as previously stated, generates an algorithm which 
solves an acceptable 
number 
of instances of the problem. 
The term ``acceptable'' in this instance is taken to mean a higher number of instances than the original EA of \cite{craven2016parallel} 
with $H_2$ inserted (cf. `$H_2$' column of 
Tables \ref{basic_tests}--\ref{basic_tests_3} given later). 
To recap, our hyper-heuristic controls the injection of heuristic chains into an EA in order to determine the best heuristic chain. 
The initial heuristic chain can be the best heuristic known (i.e., $H_2$) or a random chain. 
If the initial heuristic chain is random, then the heuristic generator is called. 
This random chain is set to a random length between 2 and 10.
We now present the core algorithmic contribution and how these algorithms are related. Algorithm \ref{alg_1} tests heuristic chains; the parameters used are $H_{\max}=20$, $N_{train}=15$, $N_{test}=50$ and $N_{valid}=50$. 

\begin{algorithm}[h!]
\caption{Heuristic generation and testing methodology }\label{alg_1}
\begin{algorithmic}[1]

\Statex \textbf{Input: } Group $G$; parameters: number of training instances $N_\textrm{train}$; number of testing instances $N_\textrm{test}$; number of validation instances $N_\textrm{valid}$; initial heuristic chain; maximum number, $C_{\max}$, of heuristics to generate.
\Statex \textbf{Output: } Runtime statistics; best heuristic chain found.

\State $i^*\gets 0$, \,$i\gets 1$
\While{$i\leq C_{\max}$}
\If{$i=1$}
\State{$C_i\gets\textrm{initial chain}$}
\Else
\State{Call heuristic chain generator (Algorithm \ref{heur_gen}), giving chain $C_i$.}
\EndIf
\State{Execute the EA, with injected chain $C_i$ on all training instances.}
\Statex\hspace{1cm}\Comment{Get metric $M_{i,train}$ on training set.}
\If{$i=1$}
\State{$M_{train}*\gets M_{i,train}$, \,\,\,$i^*\gets i$}
\Else
\If{$M_{i,train}<M_{train}*$}\Comment{Better chain found for training set;}
\Statex\hspace{6.65cm}test chain on the testing set.
\State{Execute the EA, with injected chain $C_i$ on all test instances.}
\Statex
\Comment{Get metric $M_{i,test}$.}
\State{If $M_{1,test}$ does not exist then execute the EA with injected chain}
\Statex \hspace{1.7cm}$C_1$ on all testing instances. Let $M_{test}^*\gets M_{1,test}$.
\If{$M_{i,test}<M_{test}^*$}
\State{$i^*\gets i$}\Comment{A better chain has been found on the testing set.}
\EndIf
\Else
\State{Accept chain $C_i$ (i.e., $i^*\gets i$ ) with probability $p_h$.}

\Comment{Otherwise rewind chain back to the last best chain.}
\EndIf
\EndIf
\EndWhile
\If{$i^*\neq 1$}
\State{Compare chain $C_p$ with chain $C_1$ on the validation set of instances 
\Statex\hspace{0.5cm} via execution of the EA with injected chains (i) $C_p$ and (ii) $C_1$.}
\State{\textbf{return} timeout and $C_i*$. End.}
\EndIf
\end{algorithmic}
\end{algorithm}

The EA referred to in Algorithm \ref{alg_1} is the EA of \cite{craven2016parallel} run on an input collection of instances.
The EA parameter values are reduced as in Table \ref{new_EA_params}. 
Note also that there is a probability, $p_h$, that the current chain will be accepted if it does not perform better than the best chain found (on the training instances) so far.

The group definition of $G$ is a piece of code which simply defines the group, 
its instance parameters over which the instance will be computed, and the cost functions. 
Next is the heuristic chain generator, Algorithm \ref{heur_gen}.
If the initial heuristic chain is a random chain, then this random chain is created by appending a given number 
(here, a random number between two and ten) of simple heuristics randomly chosen from $H_1$--$H_7$.
Otherwise, the heuristic generator (Algorithm \ref{heur_gen}) generates new chains of simple heuristics from the chain given by the current step of Algorithm \ref{alg_1} 
by a process of insertion, deletion or substitution at random positions in the heuristic chain. 
The heuristic is then returned in the form of a series of commands written into a file read by the EA when it is time to execute the chain. 
Chains not allowed include the set of all chains of the form $H_3^k$ for some $k>0$ (i.e., a chain consisting solely of deletions) 
or chains that are identical to those already examined in the hyper-heuristic run. We let $p_i=p_s=0.4$ and $p_d=0.2$.

\begin{algorithm}[h!]
\caption{Heuristic chain generator}\label{heur_gen}
\begin{algorithmic}[1]
\Statex \textbf{Input: Set of heuristic chains $\mathcal{C}=\left\{C_1,\ldots,C_k\right\}$ already examined.}
\Statex \textbf{Output: New heuristic chain $C'$.}
\Statex $C'\gets C_i$, the heuristic chain given by Algorithm 1.
\While{$C'\in \mathcal{C}$}
\Statex Choose operation at random subject to probabilities $p_i,p_s,p_d$ (of insertion, substitution and deletion respectively).
\Statex Perform chosen operation on $C'$ with a simple heuristic chosen at random from $H_1$--$H_7$ (if not deletion).
\EndWhile
\State{\textbf{return} heuristic chain $C'$. End.}
\end{algorithmic}
\end{algorithm}

\label{inst_gen}
An instance generator is also used. 
This creates instances at random, with random number seed based upon the computer clock. 
Included are instance parameters ($N$, $\ell$, $L_1$, $L_2$, $G$ - Section \ref{defn_problem}), a random word function, and the cost functions as in Section \ref{defn_problem}.

\subsection{Details of Implementation}
During early development of the hyper-heuristic,
issues with speed were noted. 
A number of measures were put into place to increase processing speed.
Firstly, an EA population size of 25 was used (with one slave processor being assigned to each population member).
In addition, smaller EA iteration limits than \cite{craven2016parallel} were set. 
On the training and testing instances, `maxsteps' is set to 50 for degrees 1, 2 and 3 of the polynomial $f$ defining the number field $K$ (which, of course, defines $G$), and 100 for degrees 5 and 7.
On the validation instances, `maxsteps' is set to 1250 for degrees 1, 2, and 3, and 2500 for degrees 5 and 7);
this had the effect of a small decrease in the success rate of the EA compared to that of \cite{craven2016parallel} (and so the results are not directly comparable). 
The polynomials $f$ used for the above degrees were $x-1$, $x^2-x-1$, $x^3-x-1$, $x^5-x^3-1$ and $x^7-x^3-1$, being consistent with those of \cite{garber2015length,craven2016parallel}.

All instances are run with an initial word length of 10 generators (in EA generation 1) to avoid bias to the insertion operators which would occur with an initial length of 1 (for example).
No instances of degree 9 or above were attempted due to the time complexity of computation and reduction of words in the groups concerned (for more details the interested reader should consult \cite{craven2016parallel}). 
The number of instances used in each phase of the hyper-heuristic were fifteen (training), fifty (testing) and fifty (validation). 
The number of heuristic chains run by the hyper-heuristic is $H_{\max}=20$.

All experiments were run on a high-performance cluster containing Intel Xeon E5620 CPU processors, each running at 2.40 GHz. 
The hyper-heuristic was implemented in the GAP language \cite{GAP4}, and the \verb=Polycyclic= \cite{Polycyclic} package for GAP was used for computation with polycyclic group elements. 
The ParGAP \cite{ParGAP} package was also used to handle MPI communications between processors. 
Due to the domain, the popular hyper-heuristic packages such as Hyflex \cite{hyflex2012} are not suitable for use
because we are using GAP, a specialist group theory language.
As above, each experiment was run on 26 cores (1 `master' core to control, and 25 `slave' cores, one for each EA population member). 
The code referred to in this section is available from 
\url{https://github.com/MJCraven/Hyperheuristic_group}, with the instances available at \cite{cravenwoodward2020}.

\section{Experimental Results}\label{experimental results}

In this section,
hyper-heuristic
experiments are run, varying initial input and instance parameters. 
To recap, the EA with the heuristic chains injected is then executed on the previously detailed fifteen training instances. 
If the performance improves over that of previous heuristic chains then it is run with the testing set (fifty random instances). 
If the performance over this set improves over that of previous heuristic chains then the current chain is
assigned as the new best chain.
This is continued until the end of the run,
after which the chain is validated over the validation set of (a distinct set of) fifty random instances. 
For a single hyper-heuristic run, for each of twenty heuristic chains and, assuming at least one better heuristic chain is found, around 500 problem instances are run are total.

\subsection{The Best Simple Heuristic}

A LBA attack (i.e., a hillclimber) was created for each simple heuristic. 
These attacks were run on a selection of random instances, with the percentage of successful runs as 1.7\%, 51.7\%, 0\%, 0\%, 1.7\% and 1.7\% respectively for $H_1$-$H_6$. This indicates that a heuristic on its own, unless it builds appropriate solutions, is unlikely to be successful for a large set of random instances. 
In this case, $H_2$ seems to be more successful since it builds solutions by gradually increasing solution length. 
Hence, the hyper-heuristic is initialised with the chain composed solely of a single execution of $H_2$.

\subsection{Observations on the Evolution to Build Heuristic Chains}

The following details are presented for each experiment. 
The first column of Tables \ref{basic_tests}--\ref{basic_tests_3} is the degree of the polynomial $f$,
one of the main instance parameters.
The second column is the validation set metric 
(success rate, mean cost from unsuccessful runs, mean number of generations from successful runs to solve the instance) from the EA with the best known heuristic chain (insertion - $H_2$).
The third column is the validation set metric from the EA with best injected heuristic chain found,
followed by the iteration on which the best heuristic chain was found.
The fifth column gives the chain, where $H_i^k$ refers to $k$ repeated executions of heuristic $H_i$.
The last column is the number of hyper-heuristic runs it took to find the best heuristic; unsuccessful runs were those for which either either no better chain than $H_2$ was found, or, more commonly a better chain (on the grounds of testing and training performance) was found but performed worse than $H_2$ on the validation set.

\setlength{\tabcolsep}{2pt}

\begin{table}[h!]
\begin{center}
\begin{tabular}{cccccc}
\hline
$d$ & Insertion ($H_2$) & GA with Chain & Iter & Best Chain & \# Runs \\
\hline
1 & [100\%, 0, 7.88] & [100\%, 0, 7.62] & 6 & $H_2H_1H_4$ & 4 \\
2 & [100\%, 0, 157.08] & [100\%, 0, 96.04] & 2 & $H_2H_7$ & 1 \\
3 & [100\%, 0, 101.84] & [100\%, 0, 81.54] & 16 & $H_5H_3^2H_7H_2H_5$ & 2 \\
5 & [60\%, 299.55, 491.23] & [66\%, 329.53, 695.39] & 10 & $H_3H_7$ & 6 \\
7 & [32\%, 476.44, 785.94] & [42\%, 557.90, 854.05] & 6 & $H_6H_3H_7$ & 1 \\
\hline
\end{tabular}
\caption{Comparison of results on fifty validation instances.
The parameters used were $N=20$, $L_1=10$, $L_2=13$, $L=5$, as in Section \ref{inst_gen}. Those instances used by the present work are taken from the same distributions as those used by \cite{craven2016parallel}.}
\label{basic_tests}
\end{center}
\end{table}

\begin{table}[h!]
\begin{center}
\begin{tabular}{cccccc}
\hline
$d$ & Insertion ($H_2$) & GA with Chain & Iter & Best Chain & \# Runs \\
\hline
1 & [100\%, 0, 3.92] & [100\%, 0, 3.60] & 15 & $H_7H_1H_2$ & 12 \\
2 & [100\%, 0, 38.30] & [100\%, 0, 32.68] & 18 & $H_2H_3$ & 1 \\
3 & [100\%, 0, 80] & [100\%, 0, 66.60] & 5 & $H_2H_1H_2H_4$ & 1 \\
5 & [76\%, 25.67, 501.47] & [92\%, 20.75, 488.26] & 13 & $H_7H_3H_2H_1H_5H_4H_7H_5$ & 1 \\
7 & [58\%, 37.33, 585.34] & [66\%, 41.29, 497.39] & 7 & $H_5H_2H_7H_3^2$ & 6 \\
\hline
\end{tabular}
\caption{Comparison of results on fifty validation instances. The parameters used were $N=5$, $L_1=5$, $L_2=8$, $L=5$.}
\label{basic_tests_2}
\end{center}
\end{table}

\begin{table}[h!]
\begin{center}
\begin{tabular}{cccccc}
\hline
$d$ & Insertion ($H_2$) & GA with Chain & Iter & Best Chain & \# Runs \\
\hline
1 & [100\%, 0, 7.30] & [100\%, 0, 7.02] & 9 & $H_5^4H_4H_1$ & 3 \\
2 & [96\%, 35, 141.25] & [98\%, 29, 163.51] & 3 & $H_2H_7$ & 4 \\
3  & [92\%, 37.5, 180.54] & [96\%, 37, 160.33] & 20 & $H_3H_5^2H_3H_5H_3H_1H_5H_6$ & 1 \\
5 & [52\%, 617.17, 577.38] & [58\%, 141.33, 888.03] & 8 & $H_6H_3H_4H_1$ & 2 \\
7 & [12\%, 344.64, 947.5] & [18\%, 289.76, 1115.89] & 13 & $H_4^2H_5H_6H_7H_3H_2$ & 2 \\
\hline
\end{tabular}
\caption{Comparison of results on fifty validation instances. The parameters used were $N=5$, $L_1=15$, $L_2=18$, $L=5$.}
\label{basic_tests_3}
\end{center}
\end{table}

Some observations on the results are noted in the next subsection.

\subsection{Observations and Discussion of Results}
Demonstrated through Tables \ref{basic_tests}--\ref{basic_tests_3}, 
it is clear that the approach enables the creation of more successful heuristic chains than the EA of \cite{craven2016parallel}. 
Since the hyper-heuristic relies on a stochastic algorithm (the EA),
some runs are more successful than others. 
For example, 
some hyper-heuristic runs may uncover several chains 
proving more successful than the initial heuristic chain 
(e.g., Table \ref{basic_tests_2} with $d=7$). 
On the other hand, however, 
some hyper-heuristic runs may discover no chains at all that are more effective than the initial heuristic 
(recall this information is recorded in the final column of Tables \ref{basic_tests}--\ref{basic_tests_3}). 
This latter conclusion seems to be more common for $d=1$
where a high percentage of instances are solved by the EA with the initial heuristic $H_2$.

Note, in addition, that for many small $d$ (e.g., $d=1$ or $d=2$),
all problem instances are solved by the EA with the injected simple heuristic $H_2$. 
Thus, the only option to improve performance, 
in the sense it is measured in this work,
is to solve those instances in a smaller mean number of generations. 
For example, the degree $d=1$ on Table \ref{basic_tests}, shows that 100\% of problem instances are solved by the initial heuristic in a mean of 7.88 generations. 
This is improved marginally, solving all instances with a mean of 7.62 generations by the later chain $H_2H_1H_4$. 
This suggests that for larger $d\geq 5$, for example, more 'room for improvement' is possible by the hyper-heuristic.

As is often the case with EAs and hyper-heuristics, high performance computing is an advantage due to the large amounts of time required to solve a large number of instances. 
The parameter with the largest influence is the degree $d$ (see \cite{craven2016parallel} for further details).
All the above results exhibit an improvement over the results of \cite{craven2016parallel} (and so \cite{garber2015length}). By the above results, there do not seem to be patterns formed in the best heuristic found and so it is probable that there do not exist chains that work better for one particular degree.

\subsection{Characteristics of the Framework}

Through experience, and by the above analysis, some characteristics of the framework (in the context of the AAG problem and polycyclic groups defined by a number field) are observed.

First, to the best of the authors' knowledge, 
random instances have not been classified in terms of difficulty. 
For example, 
an EA that solved instance A of a problem in an average of 3000 generations (over, say, ten repetitions) 
may well solve instance B,
with the same instance parameters, in 100 generations. 
That is, for a given set of instance parameters ($N$, $L$, $L_1$ and $L_2$) there is a large variability in difficulty for randomly-generated problems.
In the experience of the authors, this effect seems to worsen for higher degrees.
Recall that $L$ is the key length in the subgroup $\mathcal{A}\leq G$. 
Due to the lengths $L_1$ and $L_2$ (in $G$) of elements in $\mathcal{A}\leq G$ the length of the key may be rather large after mapping to its image in $G$.
Combined with the relator lengths in the presentations of the groups, 
this makes problem hardness difficult to classify.
This imposes a constraint on the hyper-heuristic, since a consistent measure of performance over a small number of instances is difficult to obtain. 
Hence, a relatively large number of instances are needed, 
at least on the testing and validation.

Combinatorial optimisation problems typically have an objective function 
where, when a small change is made to the input,
there is a correspondingly small change in the output value. 
This is reflected in the so-called ``deep-valley hypothesis'' 
\cite{doi:10.1057/jors.2010.116}.
This property is often assumed when metaheuristics are applied, as metaheuristics typically make a small change 
to the solution in order to bring about a small improvement in the objective values. 
However, the objective function in this paper is, because of the group presentations used, unlikely to display the deep-valley property
and this means that the feedback provided by a more rugged ``landscape''
does not guide the search as efficiently.
This is manifested by a heuristic chain having a low success rate on the training instances but also having a high success rate on the testing instances, or vice-versa. 

Hyper-heuristics may be applied to continuous optimisation problems, 
where real-valued feedback from the objective function may guide the search process. 
The situation is more complex for the current optimiser since the objective value is discrete: 
that is, the optimiser has either solved a given instance or it has not. 
This work goes some way to ameliorate this issue by including the least EA cost reached as part of the performance metric. 
The hyper-heuristic is hill-climbing in the space of heuristic chains. 
In the next section, the paper is concluded.

\section{Concluding remarks}
\label{Concluding remarks}

This work exhibits the automatic generation of novel heuristic chains to improve
an existing EA which has previously been demonstrated to effectively attack a given KEP. 
That is, this approach is a framework for learning
(i.e., generating and testing in a hyper-heuristics setting)
cryptanalytic attacks. 
We are not proposing a \emph{single} algorithm
to tackle this problem as many previous papers have done. 
Our stance is distinctly different: 
we propose a \emph{framework} to automatically generate 
algorithms for the attack. 
One of the advantages of this approach is that 
it automates the rather mechanical task of 
generating new attack algorithms for which there are often few design principles to guide us. 
This is thus an ideal match for a 
generative hyper-heuristics approach
where novel algorithms can be freely and easily generated. 
This avoids the task of manually
generating algorithms for which we often have scant
means of evaluating their effectiveness
other than actually testing them out on problems of interest. 
An evaluation metric is all a hyper-heuristic need to produce new heuristic chains. 
This article makes the following key contributions to the field:

\begin{enumerate}
\item
the proposal that hyper-heuristics is a suitable framework
in which to generate and test heuristic chains to break a KEP. 
\item the implementation and application of a hyper-heuristic to automatically build chains of simple heuristics to break a given KEP. 
\item the demonstration that chains of simple heuristics trained on one set of problem instances can then generalise to solve a second independent set of problem instances. 
\end{enumerate}

In the realms of further work, 
we would like to generalise the hyper-heuristic framework to work towards proving or disproving security of proposed group-theoretic KEPs.
The framework exhibited is expandable, 
enabling other groups and group-theoretic problems to be used.
Possible other uses could be to show that some proposed KEPs prove resistant to LBA attacks 
(that is, KEPs for which the hyper-heuristic does not yield high-performing heuristic chains after many runs). 
For example, would the conjugacy search problem in finitely generated metabelian groups or generalised metabelian Baumslag-Solitar groups \cite{gryak2019conjugacy} 
be breakable (the authors in a preprint suggest that LBA algorithms are ineffective) with the approach? 
Similarly, would the $n$-root and subgroup membership search problems in polycyclic groups \cite{gryak2016status}, or the conjugacy problem and hidden subgroup problem in Engel groups \cite{kahrobaei2020algorithmic}, be breakable?
Further, by \cite{doi:10.1080/10586458.2018.1434704}, there are
open questions related to complexity of some problems in polycyclic groups 
(power conjugacy problem, geodesic length problem, $n$-root problem, or the subgroup membership search problem) or other problems that may be used in KEPs.
The complexity of the above problems may be analysed using the hyper-heuristic framework, potentially giving further information about the exact solutions of these problems (if they exist).

A final area is of
determining the effectiveness of using the hyper-heuristic for 
very high parameter settings across general group-theoretic problems. In this sense, 
the hyper-heuristic approach would prove very slow (since, by the experience of the authors, the vast majority of runtime tends to spent evaluating the EA cost function) and thus research into surrogate cost function approaches would be of interest
\cite{brownlee2015metaheuristic}.

It is hoped this work may encourage machine learning and hyper-heuristic approaches in cryptology. 
This may have an impact upon post-quantum cryptography, 
with such problems as the hidden subgroup problem ripe for the attack \cite{horan2018hidden}. 
If the approach is effective, 
then it confirms a given problem combined with a platform group is breakable, 
whereas if it is not effective then this may provide further evidence to validate as ``quantum-safe'' proposed cryptographic structures.

\section*{Acknowledgments}
  \noindent The authors gratefully acknowledge the Centre for Mathematical Sciences at the University of Plymouth and the Operational Research Group at Queen Mary University of London for their generous research support and encouragement.

\bibliographystyle{plain}
\bibliography{refs}

\end{document}